\documentclass[reqno,11pt]{amsart}
\usepackage{amsmath}
\usepackage{graphicx, color, enumerate}
\usepackage{amscd}
\usepackage{amsmath}
\usepackage{amsfonts}
\usepackage{amssymb}
\usepackage{mathrsfs}
\usepackage{latexsym}

\newcommand{\be}{\begin{equation}}
\newcommand{\ee}{\end{equation}}
\newenvironment{pf}{\noindent{\it
Proof}.\enspace}{\rule{2mm}{2mm}\medskip}

\newcommand{\R}{\mathbb{R}}
 \newcommand{\Rn}{\mathbb{R}^n}

\newcommand{\E}{\mathbb{E}}

\newcommand{\h}{\mathbb{H}}
\newcommand{\dyle}{\displaystyle}
\renewcommand{\a }{\alpha }
\renewcommand{\b }{\beta }

\newcommand{\e }{\varepsilon }

\renewcommand{\l }{\lambda }
\renewcommand{\L }{\Lambda }

\newcommand{\n }{\nabla }

\renewcommand{\o }{\omega }

\newcommand{\cN}{{\mathcal{N}}}

\newcommand{\intn}{\int_{\Rn}}
\newcommand{\intR}{\int_\R}

\newcommand{\wt}{\widetilde}

\newcommand{\bu}{{\bf u}}
\newcommand{\bv}{{\bf v}}

\newcommand{\bo}{{\bf 0}}

\newcommand{\bh}{{\bf h}}

\newtheorem{Theorem}{Theorem}

\newtheorem{Lemma}[Theorem]{Lemma}

\newtheorem{remark}[Theorem]{Remark}
\newtheorem{remarks}[Theorem]{Remarks}
\newtheorem{example}[Theorem]{Example}
\newtheorem{examples}[Theorem]{Examples}

        \usepackage{latexsym}
        \usepackage{amsmath}
        \usepackage{amsfonts}
        \usepackage{amssymb}
\author[Eduardo Colorado]{Eduardo Colorado$^{*}$}
\address{\noindent Departamento de Matem\'aticas, Universidad Carlos III de Madrid,
Avenida de la Universidad 30, 28911 Legan\'es, Madrid, Spain, And
Instituto de Ciencias Matem\'aticas, ICMAT (CSIC-UAM-UC3M-UCM),
 C/Nicol\'as Cabrera 15, 28049 Madrid, Spain.}
\email{eduardo.colorado@uc3m.es, eduardo.colorado@icmat.es}
\thanks{$^*$Partially supported by Ministry of Economy and Competitiveness of
Spain and FEDER, under research project MTM2013-44123-P}

\title[Bound and ground states for some coupled NLS-KdV equations]{Existence of Bound and Ground States for a System of Coupled Nonlinear Schr\"odinger-KdV Equations}

\begin{document}
{\sf \maketitle

\noindent {\small  {\bf Abstract.}
We prove the existence of  bound and ground states for a system of coupled nonlinear Schr\"odinger-Korteweg-de Vries  equations, depending on the size of the coupling coefficient.}

\

\section{Introduction}
The aim of this note  is to show some
existence  of solutions for  a system of coupled  nonlinear
Schr\"odinger-KdV equations as follows,
\begin{equation}\label{NLS-KdV}
\left\{\begin{array}{rcl}
if_t + f_{xx} +\a fg+ |f|^2f& = & 0\\
g_t +g_{xxx} +gg_x  +\frac 12\a (|f|^2)_x & = & 0,
\end{array}\right.
\end{equation}
where $f=f(x,t)\in \mathbb{C}$ while $g=g(x,t)\in \mathbb{R}$, and $\a<0$ is the real coupling constant. System \eqref{NLS-KdV} appears in phenomena of interactions between short and long dispersive waves, arising in fluid mechanics, such as  the interactions of capillary - gravity water waves. Indeed, $f$ represents the short-wave, while $g$ stands for the long-wave; see for instance \cite{fo} and references therein.

If we define $f(x,t)=e^{i(\o t+kx)}u(x-ct)$, $g(x,t)=v(x-ct)$, with $u,\, v\ge 0$  real functions, choosing  $\l_1=k^2+\o$, $\l_2=c=2k$ and $\b=-\a$, we get that $u,\, v$ solve the following system
\begin{equation}\label{NLS-KdV2}
\left\{\begin{array}{rcl}
-u'' +\l_1 u & = & u^3+\beta uv \\
-v'' +\l_2 v & = & \frac 12 v^2+\frac 12\beta u^2.
\end{array}\right.
\end{equation}

We deal with the general case, $\l_1$ not necessarily equals to $\l_2$. We demonstrate  the existence of:

\noindent -bound states when the coupling parameter is small,

\noindent  -ground states provided the coupling factor is large, not proved before for none range of $\l_j>0$, $j=1,2$.

In the particular case  $\l_1=\l_2$ and $\b>\frac 12$ studied by \cite{dfo}, the authors proved the existence of  bound states. As a consequence of our existence results, we show that, in that range of parameters,  there exist not only bound states, if not ground states. Also, we want to point out that our method, inspired in \cite{ac1,ac2},  is different from the one in \cite{dfo}, and it seems to be more appropriate to study system \eqref{NLS-KdV2}; see Remarks \ref{rem:1}, \ref{rem:2}.

\

We use the following notation: $E$ denotes the Sobolev space $W^{1,2}(\R)$, that can be defined as the completion of $\mathcal{C}_0^1(\R)$ endowed with the norm $\|u\| =\sqrt{(u\mid u)},$
which comes from the scalar product $(u\mid w)=\intR (u'w'+uw)\,dx.$
We  denote the following equivalent norms and scalar products in $E$,
$$
\|u\|_j=\left(\intR (| u'|^2+\l_j u^2)\, dx\right)^{\frac 12},\quad ( u|v)_j=\intR ( u'\cdot v'+ \l_j uv)\, dx; \quad j=1,2.
$$
We define the product Sobolev space $\E=E\times E$. The elements in $\E$ are denoted by $\bu =(u,v)$, and $\bo=(0,0)$.
We take $ \|\bu\|=\sqrt{\|u\|_1^2+\|v\|_2^2}$ as a norm in $\E$. For $\bu\in \E$,  $\bu\geq \bo$,  $\bu>\bo$, means
that $u,v\geq 0$,  $u,v>0$ respectively. We denote $H$ as the space of even (radial) functions in $E$, and $\h=H\times H$.  We define the functional
$$
\Phi (\bu)= I_1(u)+I_2(v)- \tfrac 12\b \intR u^2v\,dx,\qquad \bu\in \E,
$$
where
$$
I_1(u)=\tfrac 12 \|u\|_1^2 -\tfrac 14\, \intR u^4dx,\qquad I_2(v)=\tfrac 12 \|v\|_2^2 -\tfrac 16\, \intR v^3dx,\qquad u,\, v\in E.
$$
We say that $\bu\in \E$ is a non-trivial {\it bound state} of \eqref{NLS-KdV2} if $\bu$ is a non-trivial  critical point
of $\Phi$. A bound state $\wt{\bu}$ is called {\it ground state} if its energy is minimal among all the non-trivial bound states, namely
\begin{equation}\label{eq:gr}
\Phi(\wt{\bu})=\min\{\Phi(\bu): \bu\in \E\setminus\{\bo\},\; \Phi'(\bu)=0\}.
\end{equation}
An expanded version of this note, with more details and further results will appear in \cite{c2}.

\section{Existence of ground states}

Concerning the ground state solutions of \eqref{NLS-KdV2}, the main result is the following.
\begin{Theorem}\label{th:1}
 There exists a real constant $\L>0$ such that for any $\b>\L$, System \eqref{NLS-KdV2} has a positive even ground state $\wt{\bu}=(\wt{u},\wt{v})$.
\end{Theorem}

We will work in  $\h$.
Setting,
$$
\Psi(\bu)=(\n\Phi(\bu)|\bu)=(I_1'(u)|u)+(I_2'(v)|v)-\frac 32\b \intR u^2vdx,$$
we define the corresponding Nehari manifold
$$
\cN =\{ \bu\in \h\setminus\{\bo\}: \Psi (\bu)=0\}.
$$
One has that
\begin{equation}\label{eq:gamma}
(\n \Psi(\bu) \mid \bu)= - \|\bu \|^2-\intR u^4<0,\quad\forall\, \bu\in \cN,
\end{equation}
and thus $\cN$ is a smooth manifold locally near any point $\bu\not= \bo$ with $\Psi(\bu)=0$. Moreover,
$\Phi''(\bo)= I_1''(0)+I_2''(0)$ is positive definite, then we infer that $\bo$ is a strict minimum for $\Phi$. As a consequence,  $\bo$ is an isolated point of the
set $\{\Psi(\bu)=0\}$, proving that  $\cN$ is a smooth complete manifold of codimension $1$, and there exists a constant $\rho>0$ so that
\be\label{eq:bound}
\|\bu\|^2>\rho,\qquad\forall\,\bu\in \cN.
\ee
Furthermore, \eqref{eq:gamma} and \eqref{eq:bound} plainly imply that $\bu\in \h\setminus\{\bo\}$ is a critical point of $\Phi$ if and only if $\bu\in\cN$ is a critical point of $\Phi$ constrained on $\cN$.

Note that by the previous arguments, the Nehari manifold $\cN$ is a natural constraint of $\Phi$. Also it is remarkable that working on the Nehari manifold,
the functional $\Phi$ takes the form:
\be\label{eq:restriction0}
\Phi|_{\cN}(\bu)= \frac 16\|\bu\|^2+\frac{1}{12}\intn u^4dx=:F(\bu),
\ee
and by \eqref{eq:bound} we have
\begin{equation}\label{eq:restriction}
\Phi|_{\cN}(\bu)\ge   \frac 16\|\bu\|^2>\frac 16 \rho.
\end{equation}
Then \eqref{eq:restriction} shows that the functional $\Phi$ is bounded from below  on $\cN$, so one can try to minimize it on the Nehari manifold $\cN$.
With respect to he Palais-Smale (PS for short) condition, we remember that in the one dimensional case, one cannot expects a compact embedding of $E$ into $L^q(\R)$ for any $2< q<\infty$. Indeed, working on $H$ (the even case) it is not true too. However, we will show that for a PS sequence we can find a subsequence for which the weak limit is a solution. This fact jointly with some properties of the Schwarz symmetrization will permit us to prove Theorem \ref{th:1}.
By the previous lack of compactness, we enunciate a measure result given in \cite{l} that we will use in the proof of Theorem \ref{th:1}.
\begin{Lemma}\label{lem:measure}
If $2<q<\infty$, there exists a constant $C>0$ such that
\begin{equation}\label{eq:measure}
\intn |u|^q \, dx\le C\left( \sup_{z\in\R}\int_{|x-z|<1}|u(x)|^2dx\right)^{\frac{q-2}{2}}
\| u\|^2_{E},\quad \forall\: u\in E.
\end{equation}
\end{Lemma}

\

Let $V$ denotes the unique positive even solution of $-v''+v=v^2$, $v\in H$.
Setting
\be\label{eq:segunda}
V_2(x)=2\l_2\,V(\sqrt{\l_2}\,x)=3\l_2\mbox{ sech}^2\left(\frac{\sqrt{\l_2}}{2}x\right),
\ee
one has that $V_2$ is the unique positive  solution of $-v''+\l_2v=\frac 12 v^2$ in $H$. Hence
$\bv_2 := (0,V_2)$ is a particular solution of  \eqref{NLS-KdV2} for any $\b\in\R$.
We also put
$$
\cN_2 =\left\{v\in H : (I_2'(v)|v)=0\right\}=\left\{v\in H : \|v\|_2^2 -\frac 12\intR v^3dx=0\right\}.
$$
Let us denote $T_{\bv_2}\cN$ the tangent space of $\bv_2$ on $\cN$. Since
$$
\bh=(h_1,h_2)\in  T_{\bv_2}\cN \Longleftrightarrow
(V_2|h_2)_2=\frac 34\intR V_2^2h_2\,dx,
$$
it follows that
\begin{equation}\label{eq:tang1}
(h_1,h_2)\in T_{\bv_2} \cN  \Longleftrightarrow h_2\in T_{V_2} \cN_2.
\end{equation}
\begin{Lemma}\label{lem:gl3}
There exists $\L>0$ such that for  $\b>\L$, then $\bv_2$ is a  saddle point of $\Phi$ constrained on $\cN$.
\end{Lemma}
\begin{pf}
One has that for $\bh\in  T_{\bv_2}\cN$,
\be\label{eq:segundaa}
\Phi''(\bv_2)[\bh]^2 =\|h_1\|_1^2 +I''(V_2)[h_2]^2-\b\intR V_2 h_1^2.
\ee
According to \eqref{eq:tang1}, $\bh=(h_1,0)\in T_{\bv_2}\cN$ for any $h_1\in H$. Defining
\be\label{eq:Lambda}
\L=\inf_{\varphi\in H\setminus\{ 0\}}\frac{\|\varphi\|_1^2}{\intR V_2\varphi^2},
\ee
we have that, for $\b>\L$, there exists $\wt{h}\in H$ with
$$
\L< \frac{\|\wt{h}\|_1^2}{\intR V_2\wt{h}^2}<\b,
$$
thus, taking $\bh_0=(\wt{h},0)$ in \eqref{eq:segundaa} we find
$$
\Phi''(\bv_2)[\bh_0]^2 =\|\wt{h}\|_1^2 -\b\intR V_2 \wt{h}^2<0,
$$
finishing the proof. \end{pf}
\begin{remark}\label{rem:1}
If one consider $\l_1=\l_2$ as in \cite{dfo}, taking $\bh_0=(V_2,0)\in T_{\bv_2}\cN$ in the proof of Lemma \ref{lem:gl3} one finds that $$\Phi''(\bv_2)[\bh_0]^2 =\| V_2\|_2^2 -\b\intR V_2^3dx=(1-2\b)\| V_2\|_2^2<0\quad \mbox{provided}\quad \b>\frac 12.$$ See also Remark \ref{rem:2}.
\end{remark}
{\it Proof of Theorem \ref{th:1}.}
We start  proving that $\inf_{\cN}\Phi$ is achieved at some positive function $\wt{\bu}\in\h$. To do so, by the Ekeland's variational principle in \cite{eke}, there exists a PS sequence $\{\bu_k\}_{k\in \mathbb{N}}\subset\cN$, i.e.,
\be\label{eq:PS1}
\Phi (\bu_k)\to c=\inf_{\cN}\Phi,\qquad\n_{\cN}\Phi(\bu_k)\to 0.
\ee
By \eqref{eq:restriction0} one finds that $\{\bv_k\}$ is a bounded sequence on $\h$, and without relabeling, we can assume that
$\bu_k\rightharpoonup \bu$ weakly in $\h$, $\bu_k\to \bu$ strongly in $\mathbb{L}^q_{loc}(\R)=L^q_{loc}(\R)\times  L^q_{loc}(\R)$ for every $1\le q<\infty$
and $\bu_k\to \bu$ a.e. in $\R^2$.
Moreover,  the constrained gradient  $\n_{\cN}\Phi (\bu_k)=\Phi' (\bu_k)-\eta_k \Psi'(\bu_k)\to 0$,  where $\eta_k$ is the corresponding Lagrange multiplier. Taking the scalar product with $\bu_k$ and recalling that  $(\Phi'(\bu_k)\mid \bu_k)=\Psi(\bu_k)=0$, we find that
$\eta_k (\Psi'(\bu_k)\mid \bu_k)\to 0$ and this
jointly with \eqref{eq:gamma}-\eqref{eq:bound} imply  that $\eta_k\to 0$. Since in addition $\|\Psi'(\bu_k)\|\leq C<+\infty$, we deduce that
$\Phi'(\bu_k)\to 0$.

Let us define $\mu_k=u_k^2+v_k^2$, where $\bu_k=(u_k,v_k)$.
By Lemma \ref{lem:measure}, applied in a similar way as in \cite{c-fract}, we can prove that there exist $R, C>0$ so that
\be\label{eq:vanishing}
\sup_{z\in\R}\int_{|z|<R}\mu_k\ge C>0,\quad\forall k\in\mathbb{N}.
\ee
We observe that   we can find a sequence of points
$\{z_k\}\subset\R^2$ so that by \eqref{eq:vanishing}, the translated sequence $\overline{\mu}_k(x)= \mu_k(x+z_k)$ satisfies
$$
\liminf_{k\to\infty}\int_{B_R(0)}\overline{\mu}_k\ge C >0.
$$
Taking into account that $\overline{\mu}_k\to \overline{\mu}$ strongly in $L_{loc}^1(\R)$,
we obtain that $\overline{\mu}\not\equiv 0$. Therefore, defining
$\overline{\bu}_k(x)=\bu_k(x+z_k)$, we have that  $\overline{\bu}_k$ is also
a PS sequence for $\Phi$ on $\cN$, in particular the weak limit of $\overline{\bu}_k$, denoted by $\overline{\bu}$,
is a non-trivial critical point of $\Phi$ constrained on $\cN$, so $\overline{\bu}\in\cN$. Thus, using \eqref{eq:restriction0} again, we find
$$
\Phi (\overline{\bu}) = \dyle F(\overline{\bu}) \le \dyle\liminf_{k\to\infty} F(\overline{\bu}_k)= \dyle\liminf_{k\to\infty}\Phi(\overline{\bu}_k)= c.
$$
Furthermore, by Lemma \ref{lem:gl3} we know that necessarily
$\Phi(\overline{\bu})<\Phi(\bv_2)$.
Clearly  $\wt{\bu}=|\overline{\bu}|=(|\overline{u}|,|\overline{v}|)\in\cN$ with
\be\label{eq:min}
\Phi(\wt{\bu})=\Phi (\overline{\bu})=\min\{\Phi(\bu)\, :\: \bu\in\cN\},
\ee
and $\wt{\bu}\ge \bo$. Finally, by the maximum principle applied to each single equation and the fact that $\Phi(\wt{\bu})<\Phi(\bv_2)$, we get $\wt{\bu}> \bo$.

\

To finish, one  can use  the classical properties of the Schwartz symmetrization to each component,  proving that  $\wt{\bu}$ is indeed a ground state of \eqref{NLS-KdV2}, i.e.,
\begin{equation}\label{eq:sym}
\Phi(\wt{\bu})=\min\{ \Phi(\bu)\, :\: \bu\in\E,\:\Phi'(\bu)=0\}.
\end{equation}
\rule{2mm}{2mm}\medskip
\begin{remark}\label{rem:2}
As we anticipate at the introduction, see also Remark \ref{rem:1}, in the range of parameters by \cite{dfo}, $\l_1=\l_2$ and $\b>\frac 12$, we have found ground state solutions in contrast with the bound states founded by \cite{dfo}.
\end{remark}

\vspace{2mm}

\

\section{A perturbation result. Existence of Bound states}

\

Finally, we establish  existence of bound states to \eqref{NLS-KdV2}, provided the coupling parameter is small.
Let us set $\bu_0=(U_1,V_2)$, where $V_2$ is given by \eqref{eq:segunda} and $U_1(x)=\sqrt{2\l_1}\mbox{ sech}\, \left(\sqrt{\l_1}x\right)$ is the unique positive solution of $-u''+\l_1u=u^3$ in $H$. Then we have the following.

\begin{Theorem}\label{th:2}
There exists $\e_0>0$ such that for any $0<\e <\e_0$ and $\b=\e\wt{\b}>0$, System \eqref{NLS-KdV2} has an even bound state $\bu_{\e}> \bo$ with  $\bu_{\e}\to \bu_0$ as $\e\to 0$.
\end{Theorem}
In order to prove this result, we can follow some ideas of the proof of \cite[Theorem 4.2]{c} with appropriate modifications.
To be short, the idea is that by the non-degeneracy of $U_1$ and $V_2$ as critical points of their corresponding energy functionals on the radial space $H$,  plainly $\bu_0$ is a non-degenerate critical point of $\Phi$ on $\h$, hence, an application of the local inversion theorem and some energy computations permit us to prove the existence of $\e_0>0$ and a convergent sequence $\bu_\e\to\bu_0$ as $\e\to 0$ for $0<\e<\e_0$. It remains to show the positivity of $\bu_\e$ which relies on variational techniques in a similar way as in \cite{c}, with appropriate changes.

\end{document}